\documentclass[11pt]{article}
\usepackage{amsmath,amsfonts,graphicx,psfrag,fancybox,qsymbols,indentfirst}
\usepackage[latin1]{inputenc}
\def \be{\begin{eqnarray*}}
\def \ee{\end{eqnarray*}}
\def \ben{\begin{eqnarray}}
\def \een{\end{eqnarray}}

\def\AArm{\fam0 }
\def\AAk#1#2{\setbox\AAbo=\hbox{#2}\AAdi=\wd\AAbo\kern#1\AAdi{}}%
\def\AAr#1#2#3{\setbox\AAbo=\hbox{#2}\AAdi=\ht\AAbo\raise#1\AAdi\hbox{#3}
}%
\def \eref#1{(\ref{#1})}
\def \sur#1#2{\mathrel{\mathop{\kern 0pt#1}\limits^{#2}}}
\def\BBc{{\AArm C\AAk{-1.02}{C}\AAr{.9}{I}{\AAFf\char"3F}}}%
\def\BBn{{\AArm I\!N}}%
\def\BBp{{\AArm I\!P}}%

\def\BBr{{\AArm I\!R}}%
\def\BBt{{\AArm T\AAk{-.62}{T}T}}%
\def\BBz{{\AArm Z\!\!Z}}%
\def\BBone{{\AArm 1\AAk{-.8}{I}I}}%

\def\E{{\math E}}

\def\h{\frac{1}{2} \BBz^-}

\def\Brp{\Big)}
\def\Blp{\Big(}
\def\brp{\big)}
\def\blp{\big(}

\def\E{\ensuremath{\mathbb{E}}}

\def\videbox{\mathbin{\vbox{\hrule\hbox{\vrule height1ex \kern.5em\vrule height1ex}\hrule}}}

\def \build#1#2#3{\mathrel{\mathop{\kern 0pt#1}\limits_{#2}^{#3}}}

\def\proof{\noindent{\bf Proof:}\hskip10pt}
\def\QED{\hfill\vrule height 1.5ex width 1.4ex depth -.1ex \vskip20pt}

\font\cls=callig15 scaled 1000
\DeclareMathOperator{\und}{UNDER}

\def\cro#1{\llbracket #1  \rrbracket}

\begin{document}
\newtheorem{lem}{Lemma}[section]
\newtheorem{defi}[lem]{Definition}
\newtheorem{theo}[lem]{Theorem}
\newtheorem{cor}[lem]{Corollary}
\newtheorem{prop}[lem]{Proposition}
\newtheorem{nota}[lem]{Notation}
\newtheorem{rem}[lem]{Remark\rm}
\newqsymbol{`P}{\mathbb{P}}
\newqsymbol{`E}{\mathbb{E}}
\newqsymbol{`N}{\mathbb{N}}
\newqsymbol{`R}{\mathbb{R}}
\newqsymbol{`o}{\omega}
\newqsymbol{`M}{{\cal M}}

\newdimen\AAdi%
\newbox\AAbo%
\font\AAFf=cmex10

\setcounter{page}{1}
\setcounter{section}{0}\date{}

\begin{center}
\LARGE{\bf Martingales and Profile of Binary  Search Trees }\\
\medskip\medskip{\Large B. Chauvin, T. Klein, J-F. Marckert , A. Rouault,}\\
\medskip\medskip\medskip
\normalsize
Universit\'e de Versailles \\
45 Avenue des Etats Unis\\
78035 Versailles Cedex\\
France
\end{center}
\medskip\medskip

\normalsize\rm



\begin{abstract}
We are interested in  the asymptotic analysis of the binary search tree (BST) under the random permutation model. 
Via an embedding  in a continuous time model, 
 we  get new results, in particular 
the asymptotic behavior of the profile.
\end{abstract}

 \vskip 5mm
{\bf A.M.S. Classification.}  {\tt 60J25, 60J80, 68W40, 60C05, 60G42, 60G44. }

\section{Introduction}\label{intro}

This paper deals mainly with two classical models of binary trees processes: 
the binary search tree process and the Yule tree process.\\
$\bullet$ A labeled binary search tree (LBST) is a structure used in computer science 
to store totally ordered data. At time $0$ the LBST is reduced to a leaf without label. 
Each unit of time, a new item is inserted in  a leaf of the tree. 
This leaf is then replaced by an internal node 
with two leaves. We are interested in the sequence 
of underlying unlabeled
trees $({\cal T}_n)_n$ induced by this construction.
We call this sequence the binary search tree process, or BST process. 
 \\
$\bullet$ The Yule tree process $(\BBt_t)_t$ is a continuous time (unlabeled) binary tree process in which 
each leaf behaves independently from the other ones (at time $0$, the tree $\BBt_0$ is reduced to a leaf).
After an (random) exponential time, a leaf has two children. Due to the lack of memory 
of the exponential distribution, each leaf is equally likely the first one to produce children.

Under a suitably chosen random model of data (the random permutation model), 
the two models of trees are deeply related. 
In the Yule tree process, let $\tau_n$ be the
 random time  
when the $n+1$th leaf appears. Under the random permutation model the link between the two models is the following one: the process $(\BBt_{\tau_n})_n$
 has the same law as 
 $({\cal T}_n)_n$.
This allows the construction of
the BST process and the Yule tree process on the same probability space on which $({\cal T}_n)_n = (\BBt_{\tau_n})_n$. 
This embedding of the BST process into a continuous time model allows to use independence properties 
between subtrees in the Yule tree process (it is a kind of Poissonization).  
 Many functionals of the BST 
 can then be derived using known results on the Yule tree.  An interesting quantity is the profile of ${\cal T}_n$ which is the sequence $(U_k(n))_{k \geq 0}$ 
 where $U_k (n)$ is the number of leaves of  ${\cal T}_n$ at level $k$.
 Here, in \eref{mc}, 
 the martingale family $(`M_n (z))_n$ -- the Jabbour's martingale -- 
 which encodes the profile of $({\cal T}_n)_n$
 is shown to be 
 strongly related to the  martingale family $(M_t (z))_t$  
 that encodes the profile of $(\BBt_t)_t$.
  
The aim of the present paper is to revisit the study of $(`M_n (z))_n$
  using  the embedding. For $z > 0$,
we recover very quickly the behavior of the limit $`M_\infty (z)$: 
 positive 
 when $z \in (z_c ^- , z_c^+)$, 
 zero when 
 $z \notin [z_c ^- , z_c^+]$. In the critical cases  $z = z_c ^\pm$ the behavior was unknown. 
 We prove that $`M_\infty (z^\pm _c) = 0$ a.s. and get
  the convergence of the derivative. 
The limits $`M'_\infty(z)$ and $`M_\infty(z)$ satisfy a splitting formula \eref{master'} which, for $z=1$ gives the Quicksort equation (Corollary \ref{cor11}).
  Thus, the embedding method is the key tool for proving and enlarging convergence results on the BST martingale (Theorem \ref{bigjab}) and its derivative (Theorem \ref{theoderiv}).
See the companion paper \cite{ChauRou03} for complements.

 The paper is organized as follows. After the definition of the models in Section 2, 
 we explore some consequences of the embedding. 
 In particular, in 2.5.1 we exhibit  a family of uniform random variables attached to the nodes of the Yule tree.
These random variables give, for every node $u$, the limiting proportion of leaves 
issued from $u$ among those issued from its parent. A similar property holds for the embedded BST.
%
In subsection 2.5.2, the appearance of uniform variables as limiting proportion of leaves 
is explained on a  LBST model.
In Section \ref{marti}, we study the convergence, as $n \rightarrow \infty$, 
of the BST martingale ${\cal M}_n (z)$.

Thanks to this method, we are able in Section \ref{prof} to describe the asymptotic behavior 
 of the profile $U_k (n)$ when $k \simeq 2z \log n$ in the whole range $z \in (z_c^- , z_c^+)$.
   Previously, the result was known only on a sub-domain 
where the $L^2$ method works (\cite{Jab2}). 

Finally, in Section \ref{tag}, an other point of view is investigated. 
In biasing the evolution rules of the BST -- a random line of descent is distinguished, 
and the evolution of the nodes belonging to this branch is different from the other ones -- 
it appears that the behavior of the nodes on the distinguished line gives information on the whole tree. 
\section{The models}
\label{topmodel}
\subsection{Binary search trees}

For a convenient definition of trees we are going to work with, let us first define 
\[\mathbb{U} =\{\emptyset\}\cup\bigcup_{n\geq 1}\{0,1\}^n\]
the set of finite  words on the alphabet $\{0,1\}$ (with $\emptyset$ for the empty word). For $u$ and $v$ in $\mathbb{U}$, denote by $uv$ the concatenation of the word $u$ with the word $v$ (by convention we set, for any $u\in\mathbb{U}$, $\emptyset u=u$). If $v \not= \emptyset$, we say that $uv$ is a descendant of $u$ and $u$ is an ancestor of $uv$.
Moreover $u0$ (resp. $u1$) is called left (resp. right) 
  child of $u$.

A \sl complete binary tree \rm $T$ is a finite subset  of $\mathbb{U}$ such that
\[\left\{
\begin {array}{l}
\emptyset\in T\\
\textrm{ if }uv\in T \textrm{ then } u\in T \,,\\
u1 \in T \Leftrightarrow u0\in T\,.
\end {array}
\right.\]
The elements of $T$ are called \it nodes\rm, and $\emptyset$ is called the \it root \rm ; \rm $|u|$, the number of letters in $u$, is the \it depth \rm of $u$ (with $|\emptyset|=0$).
Write {\bf BinTree} for the set  of complete binary trees.

A tree $T \in$ {\bf BinTree} can be described by giving the set $\partial T$ of its \it  leaves\rm, that is, the nodes that are  in $T$ but with no descendants in $T$.
The nodes of $T\backslash\partial T$ are called \it internal \rm nodes.

We now introduce labeled binary search trees (LBST), that are widely used to store totally ordered data (the monograph of Mahmoud  \cite{Mah} gives an overview of the state of the art). 

Let $A$ be a totally ordered set of elements 
named keys and for $n \geq 1$, let $(x_1 , ... , x_n)$ be picked up without replacement from $A$.
The LBST built from these data is the complete binary tree in which each internal node is associated with a key belonging to $(x_1 , ... , x_n)$ in the following way: the first key $x_1$ is assigned to
 the root. The next key $x_2$ is assigned to the left child
 of the root if it is smaller than $x_1$, or it is
 assigned to the right child of the root if it is larger than $x_1$. 
We proceed further  inserting key by key recursively.   After the $n$ first insertion, one has a labeled binary tree in which $n$ nodes own a label: these nodes are considered as internal nodes. One adds $n+1$ (unlabeled) leaves to this structure in order to get a labeled complete binary tree with $n$ internal nodes.
 
\begin{figure}[htbp]
\begin{center}
\psfrag{a}{}\includegraphics[height=0.3cm]{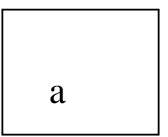}~~
\psfrag{a}{0.5}\includegraphics[height=1.2cm]{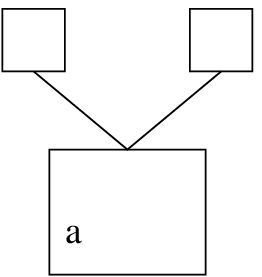}
\psfrag{c}{0.8}~~\includegraphics[height=2.1cm,width=2cm]{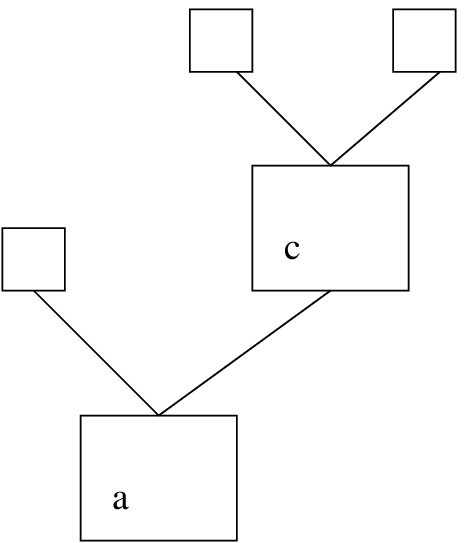}
\psfrag{g}{0.9}~~\includegraphics[width=2.6cm,height=3cm]{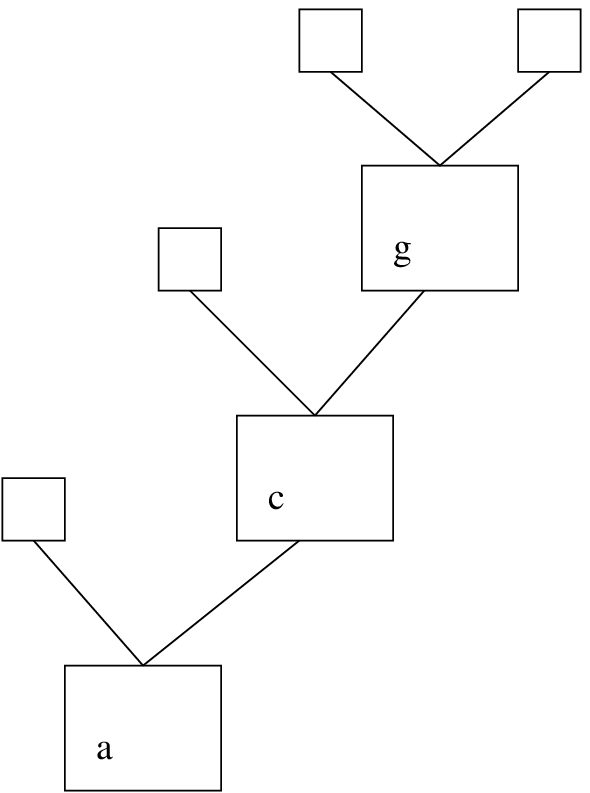}
~\psfrag{b}{0.3}\psfrag{h}{0.4}\includegraphics[height=3cm,width=3cm]{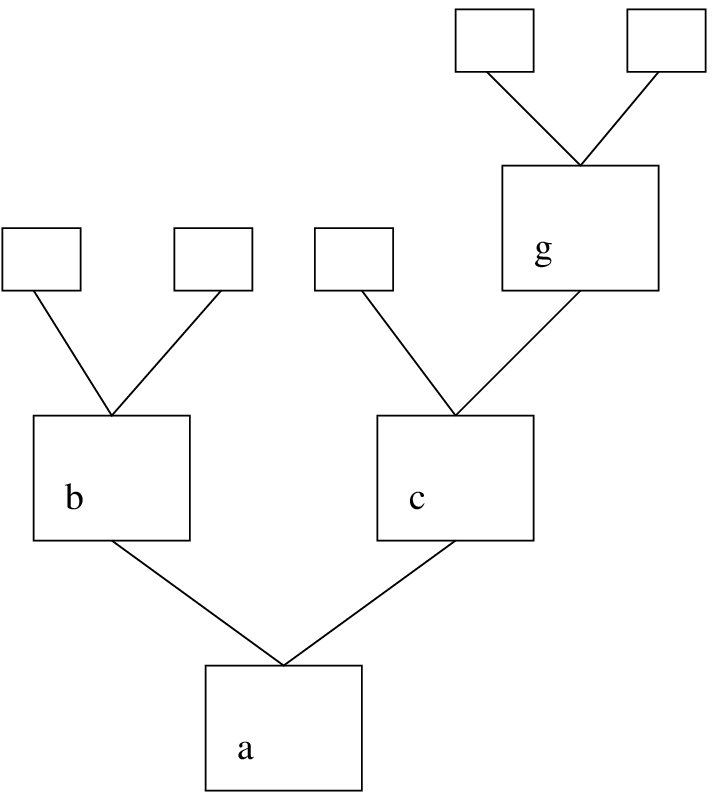}\psfrag{f}{0.4}~~\includegraphics[height=3cm,width=3cm]{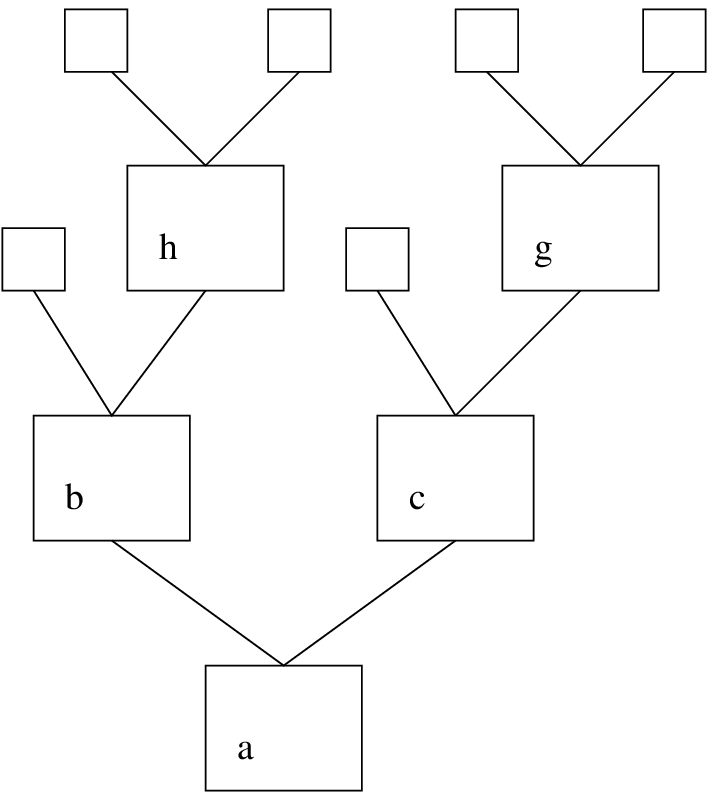}
\end{center}
\caption{BST built with the sequence of data 0.5, 0.8, 0.9, 0.3, 0.4 (empty squares are leaves).}
\end{figure}
To study the shape of these trees for large $n$, it is classical to introduce a random model.
One usually assumes that the successively inserted data $(x_i)_{i\geq 1}$ are i.i.d. random variables 
with a continuous distribution $F$. Under this model, let us call the LBST $L_n^{(F)}$; 
it is a random variable taking values in 
the set of complete binary trees
 in which each  internal node has a label in [0,1]. The sequence $(L_n^{(F)},n\geq 0)$ is a Markov chain. \par

We are, in fact, mainly interested in the underlying tree ${\cal T}_n^{(F)}$ of $L_n^{(F)}$ , i.e. the tree that 
has the same arborescent structure of $L_n^{(F)}$, but that has no label. We set
\[\big({\cal T}_n^{(F)},n\geq 0\big):=\big(\und(L_n^{(F)}),n\geq 0\big);\]
by construction  ${\cal T}_n^{(F)}$ is a complete binary tree.\par
 
For every $n \geq 1$, the string $x_1, .. ,x_n$ induces a.s. a permutation $\sigma_n$ such that 
$x_{\sigma_n(1)}<x_{\sigma_n(2)}<\dots<x_{\sigma_n(n)}$. Since the $x_i$ are exchangeable, 
 $\sigma_n$  is uniformly distributed on the set ${\cal S}_n$ of permutations of $\{1, .., n\}$. Since this claim is not sensitive to $F$ we will assume, for the sake of simplicity, that $F$ is the uniform distribution on $[0,1]$, and we write from now $L_n$ instead of $L_n^{(F)}$ and ${\cal T}_n$ instead of ${\cal T}_n^{(F)}$. 
 This is the so-called random permutation model. Again by exchangeability, 
$\sigma_n$ is independent of the vector 
$(x_{\sigma_n(1)}, x_{\sigma_n(2)}, \dots ,x_{\sigma_n(n)})$ and we have
\be
P\big(x_{n+1} \in (x_{\sigma_n(j)}, x_{\sigma_n(j+1)})\, |\,  \sigma_n \big)& =& P \big(x_{n+1} \in (x_{\sigma_n(j)}, 
x_{\sigma_n(j+1)})\big) \\
&=& P(\sigma_{n+1}(j+1) = n+1) = (n+1)^{-1}
\ee 
for every $j \in\{0, 1 , .., n\}$, where $x_{\sigma_n (0)} := 0$ and $x_{\sigma_n (n+1)} := 1$. This relation ensures the consistency of the sequence $(\sigma_n)_n$. \par
One can also express this property with the help of the sequential ranks of the permutation: the random variables 
$R_k = \sum_{j=1}^k \BBone_{x_j \leq x_k}, k \geq 1$ are independent and $R_k$ is uniform on $\{ 1, \ldots , k\}$ (see for instance Mahmoud \cite{Mah}, section 2.3), so that $P(R_{n+1} = j+1 \, |\, R_1, .. , R_n )= (n+1)^{-1}$.\par

In terms of binary search tree, this means that the insertion of the $n+1$st key in the tree with $n$ internal nodes is uniform  among its $n+1$ leaves. In other words, in the random permutation model, the sequence $({\cal T}_n)_{n\geq 0}$ is  a Markov chain on {\bf BinTree}  defined by ${\cal T}_0 = \{\emptyset\}$ and 
\ben
\label{transit}
{\cal T}_{n+1}&=&{\cal T}_n\cup \{D_{n}0,D_{n}1\}\,, \cr
P(D_{n} = u \, | \, {\cal T}_n) &=& (n+1)^{-1} , \ \ \ \ u \in \partial{\cal T}_n\,;
\een
the leaf $D_{n}$ of ${\cal T}_n$ is the random node where the $n+1$-st key is inserted, its level is $d_n$.\par

 The difference of the rule evolutions of $L_n$ (that depends deeply on the values $x_1,\dots, x_n$ already inserted) and ${\cal T}_n$ (that depends of nothing) is similar to Markov chain in random 
environment ($L_n$ is the quenched Markov chain and ${\cal T}_n$ the annealed one).

This Markov chain model  is a particular case ($\alpha = 1$) of the   diffusion-limited aggregation (DLA) on a binary tree, where a constant $\alpha$ is given and the growing of the tree is random with probability of insertion at a leaf $u$ proportional to $|u|^{-\alpha}$ (Aldous-Shields \cite{AlS}, Barlow-Pemantle-Perkins \cite{BPP}).

Here are few known results about the evolution of BST. First, the saturation level $h_n$ and the height $H_n$,
\ben
h_n = \min \{ |u| : u \in \partial {\cal T}_n \}\ \  , \ \ H_n = \max \{ |u| : u \in \partial {\cal T}_n \}
\een
grow logarithmically (see for instance 
 Devroye \cite{Dev1}
) 
\ben
\label{cc'}
\hbox{a.s.} \ \ \ 
\lim_{n \rightarrow \infty} \frac{h_n}{\log n} = c'=0.3733...\, \ \ \ 
\lim_{n \rightarrow \infty} \frac{H_n}{\log n} = c=4.31107...\,;
\een
the constants $c'$ and $c$ are the two solutions of the equation $\eta_2 (x) = 1$ where
\ben
\label{defeta}
\eta_\lambda(x) := x \log \frac{x}{\lambda} - x +\lambda, \ \ \ \ x \geq 0\,,
\een
is the Cramer transform of the Poisson distribution of parameter $\lambda$. Function $\eta_2$ reaches its minimum at $x=2$. 
It corresponds to the rate of propagation of the depth of insertion: $\frac{d_n}{2\log n} \buildrel{P}\over{\longrightarrow} 1$.
More precise asymptotics for $H_n$ can be found in \cite{Dr2}, \cite{Reed}, \cite{Robs}, \cite{KraMa}.

Detailed information on ${\cal T}_n$ is provided by the whole profile
\ben
U_k(n):= \# \{ u \in \partial {\cal T}_n , |u| = k \} \ \ , \ \ k\geq 1\,,
\een
that counts the number of leaves of ${\cal T}_n$ at each level. 
Notice that $U_k(n)= 0$ for $k >H_n$ and for $k<h_n$. To get asymptotic results, it is rather natural 
to encode the profile by the so-called polynomial level 
$\sum_k U_k (n) z^k$, whose degree is $H_n$. Jabbour \cite{Jab2,Jab1} proved a remarkable martingale property for these random polynomials. More precisely, 
for $z \notin \h = \{ 0, -1/2 , -1 , -3/2 , \cdots\}$ and $n\geq 0$,
let
\ben\label{defjab}
`M_n (z) := \frac{1}{C_n (z)} \sum_{k\geq 0} U_k (n) z^k = \frac{1}{C_n (z)} 
\sum_{u \in \partial{\cal T}_n} z^{|u|}\,,
\een
where $C_0(z)=1$ and for $n\geq 1$,
\ben
\label{defcn}
C_n (z) := \prod_{k=0}^{n-1} \frac{k + 2z}{k+1} = (-1)^n 
\begin{pmatrix} - 2z\\n\end{pmatrix}\,,
\een
and let  ${\cal F}_{(n)}$ be the  $\sigma$-field generated by 
all the events $\{ u \in {\cal T}_j\}_{j \leq n , u \in \mathbb{U}}$ .
Then $(`M_n (z), {\cal F}_{(n)})_n$ is a martingale 
to which, for the sake of simplicity, we refer from now as the 
BST 
 martingale. If $z > 0$, this positive martingale is a.s. convergent; the limit $`M_\infty (z)$ is positive a.s. if $z \in (z_c ^- , z_c ^+ )$, with
\ben
\label{defzc}
z_c ^- = c' / 2 = 0.186... , \ \ z_c ^+ = c / 2 = 2.155... 
\een 
and $`M_\infty (z) = 0$ for $z \notin [z_c ^- , z_c ^+]$ (Jabbour \cite{Jab1}). This martingale is also the main tool to prove that, properly rescaled around $2\log n$, 
the  profile has a Gaussian limiting shape (see Theorem 1 in 
\cite{Jab2}
).

\subsection{Fragmentation, Yule tree process and embedding}
 \label{YT}
 The idea of embedding discrete models (such as urn models) in continuous time branching processes goes back at least to 
Athreya-Karlin \cite{AKar}. It is described  in Athreya and Ney (\cite{ANey}, section 9) and it has been recently revisited by Janson \cite{Janson03}. For the BST, various embeddings are mentioned in Devroye \cite{Dev1}, in particular those due to Pittel \cite{Pitt86}, and Biggins \cite{JDB94,BG97}.  
Here, we work with 
a variant of the Yule process, taking into account the tree (or ``genealogical'') structure. 

First, let us define a fragmentation process $(F(t))_{t\geq 0}$ of the interval $(0,1)$ as follows: we set $I_\emptyset = (0,1)$ and for $u = u_1 u_2 ... u_k \in \mathbb{U}$, set $I_u$ the interval
$$I_u = \Big(\sum_{j=1}^k u_j 2^{-j}, 2^{-k}+\sum_{j=1}^k u_j 2^{-j} \Big).$$
Hence, each element $u$ of $\mathbb{U}$ encodes a subinterval $I_u$ of $(0,1)$ with dyadic extremities.\par
We set $F(0)=I_\emptyset = (0,1)$. An exponential $\tau_1\sim{\bf Exp}(1)$ random variable  is associated with  $F(0)$. At time $\tau_1$, the process $F.$ jumps, the interval $(0,1)$ splits in the middle into two parts and $F(\tau_1)=((0,1/2),(1/2,1))=(I_{0},I_1)$.  
After each jump time $\tau$, the fragments of $F(\tau)$ behave independently of each other. Each fragment $I_u$ splits after a ${\bf Exp}(1)$-distributed random time into two fragments: $I_{u0}$ and $I_{u1}$. Owing to the lack of memory of the exponential distribution, when $n$ fragments are present, each of them will split first equally likely. \par

\begin{figure}[htbp]
\psfrag{a}{0}
\psfrag{0}{0}\psfrag{1}{1}
\psfrag{b}{$\tau_1$}
\psfrag{c}{$\tau_2$}
\psfrag{d}{$\tau_3$}
\psfrag{e}{$\tau_4$}
\psfrag{f}{$\tau_5$}
\psfrag{g}{$t$}
\centerline{\includegraphics[height=5cm]{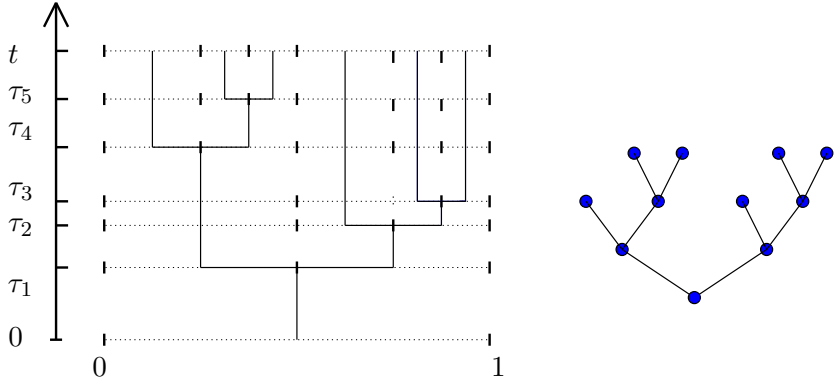}}
\caption{\label{frag2} Construction of the Yule tree from the fragmentation.}
\end{figure}

We define now the Yule tree process as an encoding of the fragmentation process. 
The idea is to interpret the two fragments $I_{u0}$ and $I_{u1}$ issued from $I_u$ 
as its two children: $I_{u0}$ is considered as the left fragment and  $I_{u1}$ the right one; like this, we obtain a binary tree structure (see Fig. \ref{frag2}).  An interval with length $2^{-k}$ corresponds to a leaf at depth $k$ in the corresponding tree structure; the size of fragment $I_u$ is $2^{-|u|}$. 
More formally, we define the tree $\BBt_t$ thanks to its set of leaves
\begin{equation}\label{fb}
\partial \BBt_t=\{u, I_u\in F(t)\}.
\end{equation}
We call Yule tree process, the \bf BinTree\rm-valued random process $(\BBt_t)_{t\geq 0}$. Both processes $(\BBt_t)_{t\geq 0}$ and  $(F(t))_{t\geq 0}$ are pure jump Markov processes.
Each process $(\BBt_t)_{t\geq 0}$ and  $(F(t))_{t\geq 0}$ can be viewed as an encoding of the other one, using (\ref{fb}) and:
\[F(t)=\{I_u, u\in \partial \BBt_t\}.\]
 The counting process $(N_t)_{t \geq 0}$ that gives the number of leaves in $\BBt_t$, 
\ben 
N_t := \# \partial\BBt_t,
\een
is the classical Yule (or binary fission) process (Athreya-Ney \cite{ANey}).

Let $0=\tau_0<\tau_1<\tau_2<...$  be the successive jump times of $\BBt.$ (or of $(F(.))$,
\begin{equation}
\label{deftau}
\tau_n = \inf \{t : N_t = n+1\}\,.
\end{equation}

The following proposition allows us to build  the  Yule tree process  and the BST on the same probability space.
This observation was also made in Aldous-Shields \cite{AlS} section 1, 
(see also Kingman \cite{JFCK} p.237 and Tavar\'e \cite{Tavar} p.164 in other contexts). 
\begin{lem}
\label{BNT}
\begin{itemize}
\item[a)] The jump time intervals $(\tau_n-\tau_{n-1})_n$ are independent and satisfy: 
\begin{equation}
\label{expo}
\tau_n-\tau_{n-1}\sim {\bf Exp}(n) \textrm{ for any } n\geq 1,
\end{equation}
where ${\bf Exp}(\lambda)$ is the exponential distribution of parameter $\lambda$.
\item[b)]
The processes $(\tau_n)_{n\geq 1}$ and $\bigl(\BBt_{\tau_n}\bigr)_{n \geq 1}$ are independent.
\item[c)] The  processes $\bigl(\BBt_{\tau_n}\bigr)_{n \geq 0}$ and $\bigl( {\cal T}_n\bigr)_{n \geq 0}$ have the same distribution.
\end{itemize}
\end{lem}
\proof  $(a)$ is a consequence of the fact that the minimum of $n$ independent random variables ${\bf Exp}(1)$-distributed is ${\bf Exp}(n)$-distributed. $(b)$ comes from the independence of jump chain and jump times. 
Since the initial states and evolution rules of the two Markov chains $\BBt_{\tau_n}$ and ${\cal T}_n$ 
are the same ones, (c) holds true. \QED

\noindent \bf Convention: (A unique probability space) \rm From now, we consider that the 
fragmentation process, the Yule tree process and the BST process are built on the same probability space. Particularly, on this space, we have
\begin{equation}\label{embe}
\bigl(\BBt_{\tau_n}\bigr)_{n \geq 0}{=} \bigl( {\cal T}_n\bigr)_{n \geq 0}.
\end{equation}
We say that the BST process is embedded in the Yule tree process.
We define the filtration $({\cal F}_t)_{t\geq 0}$ by ${\cal F}_t=\sigma(F(s),s\leq t)$. On the unique probability space, the sigma algebra ${\cal F}_{(n)}$ is equal to $\sigma(F(\tau_1),\dots,F(\tau_n))$.
\medskip

If we consider the measure valued  process $(\rho_t)_{t\geq 0}$ defined by
\ben
\label{defrho}
\rho_t=\sum_{u\in\partial \BBt_t} \delta_{-\log_2|I_u|}=\sum_{u\in\partial \BBt_t} \delta_{|u|}\,,
\een
we obtain a continuous time branching random walk. 
The set of positions is $\BBn_0 = \{ 0, 1, 2, \cdots\}$. Each individual 
 has an ${\bf Exp} (1)$ 
distributed lifetime and does not move. At his death, he disappears and is replaced by  two children, whose  positions are both their parent's position shifted by 1.
The set of  individuals alive at time $t$ is $\partial\BBt_t$ and
  the position of individual $u$ is simply $|u|$.\par
This is a particular case of the following general fact: the empirical measure of the logarithm of the size of fragments in homogeneous fragmentations with finite dislocation measures is a branching random walk (this idea goes back to Aldous and Shields \cite{AlS} Section 7f and 7g).

\subsection{Martingales and connection}

The classical\footnote{In general $|u|$ is replaced by the position $X_u$ and $L(\theta) = \beta(E  \int e^{\theta x} Z(dx) - 1)$ where $\beta$ is the parameter of lifetime and $Z$ is the offspring point process; here $\beta = 1$ and $Z = 2\delta_1$.}  family of ``additive'' martingales associated with the Yule process,  
 parameterized by $\theta$ in $\BBr$ (sometimes in $\BBc$) and indexed by $t \geq 0$, is here given by 
\[m(t, \theta) := \sum_{u \in  \partial\BBt_t
} \exp (\theta |u| - t L(\theta)) ,\]
where 
\begin{equation}\label{lth}
L(\theta) = 2e^\theta - 1
\end{equation} 
 (see \cite{Uchi}, \cite{Kyp99}, and \cite{BertRou2} for the fragmentation).
%
For easier use, we set $z = e^{\theta}$ and then 
consider the family of $({\cal F}_t,t\geq 0)$-martingales 
\begin{equation}
\label{BHM}
M(t, z) := m(t, \log z)=\sum_{u \in \partial\BBt_t} z^{|u|} e^{t (1 - 2z)}.
\end{equation}
In particular  $M(t , 1/2) = 1$ and $M(t, 1) = e^{-t} N_t$. \par
The embedding formula \eref{embe} allows to connect the family of BST martingales $({\cal M}_n,{\cal F}_{(n)})_n$ to the family of Yule martingales $(M(t,z),{\cal F}_t)_{t}$.  If we observe the martingale $(M(.,z)$ at the stopping times $(\tau_n)_n$, we can  ``extract'' (Proposition \ref{BHJ} below) the space component $`M_n (z)$ and a time component 
\ben
\textrm{\cls{C}}_n (z) := e^{\tau_n (1-2z)} C_n (z)\,.
\een
Notice that $\bigl(\textrm{\cls{C}}_n (z)\bigl)_n$ is ${\cal F}_{\tau_n}$-adapted.

A classical result (see Athreya-Ney \cite{ANey} or Devroye \cite{Dev1} 5.4)
says that, a.s., $e^{-t} N_t$ converges when $t\to +\infty$, and 
\begin{equation}
\label{limetnt}
\xi 
:= \ \lim_{t\rightarrow \infty} e^{-t} N_t\ \sim {\bf Exp}(1)\,.\end{equation}
Since $\lim_n \tau_n = \infty$ a.s. (see Lemma \ref{BNT} a) ) we get from (\ref{deftau}) and (\ref{limetnt}),
\ben
\label{etn}
\hbox{a.s.} \ \lim_n n e^{-\tau_n} = \xi\,.
\een

\begin{prop} {\bf (martingale connection)}
\label{BHJ} Let us assume $z \in \BBc \setminus\h$.
\begin{itemize} 
\item[1)]The family $\bigl(\textrm{\cls{C}}_n (z)\bigl)_{ n \geq 0}$ is a martingale with mean 1, and
\ben
\label{taun}
\hskip 1cm  \hbox{a.s.} \  \lim_n \textrm{\cls{C}}_n (z) =
 \frac{\xi^{2z-1}}{\Gamma(2z)}.
\een 
Moreover, if  $\Re z$, the real part of $z$, is positive, the convergence is in $L^1$.
\item[2)]
The two martingales $(\textrm{\cls{C}}_n (z))_{n \geq 0}$ and $(`M_n (z))_{n \geq 0}$ are independent and
\ben
\label{mc}
M  (\tau_n ,z) = \textrm{\cls{C}}_n (z)`M_n (z)\,.
\een
\end{itemize}
\end{prop}

\proof 
1) The martingale property comes from Lemma \ref{BNT} a). 
The Stirling formula gives the very useful estimate:
\ben
\label{cnz}
 C_n (z)\sim
 \frac{n^{2z-1}}{\Gamma(2z)}\,,
\een
which yields (\ref{taun}) owing to (\ref{etn}).

2) The second claim comes from \eref{embe} and \eref{BHM}, the independence comes from Lemma \ref{BNT} b).
\QED
Proposition \ref{BHJ} allows us to transfer known results about the Yule martingales to BST martingales, thus giving very simple proofs of known results 
 about the BST martingale and also getting much more.
In particular, in Theorem \ref{bigjab} 2),
we give the answer to the question asked in \cite{Jab1}, 
about critical values of $z$, with a straightforward argument.


\subsection{Limiting proportions of nodes}\label{bis} 

Let us study some meaningful random variables arising as a.s limits and 
playing an important role in 
the results of Section \ref{marti}. These variables describe the evolution of relative sizes of subtrees in Yule and BST models.

\subsubsection{In the Yule tree}

First, we recall a classical identity: let $\xi_a$ and $\xi_b$ be two independent, ${\bf Exp}(1)$-distributed 
random variables. We have
\[\frac{\xi_a}{\xi_a+\xi_b}\sim {\cal U}[0,1],\textrm{ the uniform distribution on }[0,1].\]

For every $u \in \mathbb{U}$, let $\tau^{(u)} = \inf \{ t : u \in \BBt_t \}$ 
be the time (a.s. finite) at which $u$ appears in the Yule tree, 
and for $t > 0$, set
$$\BBt^{(u)}_t =  \{ v \in \mathbb{U}  : uv \in \BBt_{t +\tau^{(u)}} \}$$
the tree process growing from $u$. 
In particular, set $$N_t ^{(u)} = \# \partial \BBt^{(u)}_t\,.$$
For $t > \tau^{(u)}$, the number of leaves at time $t$ in the subtree issued from node $u$ is 
$n_t ^{(u)} :=N_{t-\tau^{(u)}}^{(u)}$. 
The branching property and (\ref{limetnt}) give that a.s. for every $u \in \mathbb{U}$
\ben
\lim_{t \rightarrow \infty} e^{-t}
N^{(u)}_t 
= \xi_u \ \ \ ,\ \ \lim_{t \rightarrow \infty} e^{-t}
n^{(u)}_t = \xi_u \, e^{-\tau^{(u)}}\,, 
\een
where $\xi_u$ is distributed as $\xi$ i.e. ${\bf Exp}(1)$. Moreover, if $u$ and $v$ are not in the same line of descent, the r.v.  $\xi_u$ and $\xi_v$ are independent.
Since, for $t > \tau^{(u)}$, 
\ben
\label{petitn}
n_t ^{(u)} = n_t^{(u0)} + n_t^{(u1)} \ \ \ \hbox{and} \ \ \ \tau^{(u0)} = \tau^{(u1)},
\een  
a small computation yields
\ben
\label{defU}
\frac{n^{(u0)}_t}{n_t^{(u)}} \xrightarrow[]{~a.s.} U^{(u0)} := \frac{\xi_{u0}}{\xi_{u0} + \xi_{u1}} , 
\ \ \ \ \frac{n^{(u1)}_t}{n_t^{(u)}} \xrightarrow[]{~a.s.} U^{(u1)} := 1-U^{(u0)} = \frac{\xi_{u1}}{\xi_{u0} + \xi_{u1}},
\een
which allows to attach a ${\cal U}([0,1])$ r.v. to each node of $\mathbb{U}$. In particular we set
\begin{equation}\label{UUU}
U := U^{(0)} = \frac{\xi_0}{\xi_0 + \xi_1} 
\end{equation}
so that
\ben
\label{xiutau}
\xi :=\xi_{\emptyset}= e^{-\tau_1}(\xi_0 + \xi_1)\ \ ,\ \ \xi_0 = U\xi e^{\tau_1} \ \ ,\ \ \xi_1 = (1-U)\xi e^{\tau_1}\,.
\een
If $u0$ and $u1$ are brother nodes, we have $U^{(u1)} + U^{(u0)} =1$. 
We claim that
if a finite set of nodes $v_1, \ldots , v_k$ does not contain any pair of brothers, the corresponding r.v. $U^{(v_1)}, \ldots , U^{(v_k)}$ are independent. 
When none of the $v_j$ is an ancestor of another (``stopping line'' property) it is a consequence of the branching property. 
In the general case, it is sufficient to prove that $U^{(u)}$ is independent of $( U^{(v)} , v < u)$. To simplify the reading, let us give the details only for $|u| = 2$, 
for instance $u = 00$. We have, from (\ref{petitn})
$$U^{(00)} = \frac{\xi_{00}}{\xi_{00} + \xi_{01}} \ \ , \ \ 
U^{(0)} = \frac{(\xi_{00} + \xi_{01})e^{-\tau^{(00)} + \tau^{(0)}}}
{(\xi_{00} + \xi_{01})e^{-\tau^{(00)} + \tau^{(0)}} + (\xi_{10} + \xi_{11})e^{-\tau^{(10)} + \tau^{(1)}}}\,.$$
Actually, from the branching property, $\xi_{00}$ and $\xi_{01}$ are independent of $\xi_{10}, \xi_{11}, \tau^{(00)}, \tau^{(0)}, \tau^{(10)}, \tau^{(1)}$. Moreover since $\xi_{00}$ and $\xi_{01}$ are independent and  ${\bf Exp}(1)$ distributed, then $\xi_{00}/(\xi_{00} + \xi_{01})$ and 
$(\xi_{00} + \xi_{01})$ are independent, which allows to conclude that $U^{(00)}$ and $U^{(0)}$ are independent.
 
Finally, multiplying along the line of the ancestors of a node $u$, we get the representation
\ben\label{asunif}
a.s. \lim_{t \rightarrow \infty} \frac{n_t ^{(u)}}{N_t} = \prod_{v <u}  U^{(v)} \,,
\een
where the random variables $(U^{(v)})_{v\in \mathbb{U}}$ satisfy the claim.

This is of course related to multiplicative cascade models.  See \cite{ChauRou03}.

\subsubsection{In the BST}
It is straightforward to see that, by embedding, the property of the above subsection holds true for limiting proportions of nodes in the BST, as $n \rightarrow \infty$.  

This property is  also known in the LBST (this gives another proof of this property in the BST). Let us now sketch the argument for LBST.  

Assume  $x_1$ fixed. Consider the tree $L_n$ after insertion of the $n$ data $x_2,\dots, x_{n+1}$. Let $K(n):=\#\{i,i\in\cro{2,n+1}, x_i\leq x_1,  \}$ be the number of nodes in the left subtree rooted in $\emptyset$. Since the $x_i$ are i.i.d., ${\cal U}[0,1]$, the conditional distribution of $K(n)$ on $x_1$, is a binomial $B(n, x_1)$.  Hence, by the strong law of large numbers,
\[\frac{K(n)}{n}\xrightarrow[n]{a.s.} x_1.\]
Now, the subtree $t_0$ rooted in $u=0$ and the subtree $t_1$ rooted in $u=1$ are ``copies'' of $L_n$. The subtree $t_0$ is built with the random variables in the list $x_2, x_3,\dots,$ that are smaller than $x_1$ ($t_1$ is build with the ones that are larger than $x_1$). In particular, the label $x_{t_0}$ of the root of $t_0$ is the first value among $x_2, x_3,\dots,$ smaller than $x_1$. It is easy to check that $x_{t_0}$ is uniform on $[0,x_1]$, therefore it has the following representation: $x_{t_0}=x_1U$ where $U$ is uniform on $[0,1]$ and does not depend on the value $x_1$. Hence, the asymptotic proportion of nodes in the subtree $t_{00}$ is $x_1U$ while it is $x_1(1-U)$ in $t_{01}$ (what happens in the subtree $t_1$ is totally independent).\par
This iterative construction of the LBST explains why it enjoys the same property as \eref{asunif}  in the Yule process, and so does the sequence of underlying BST ${\cal T}_n$. This is a strong, which means a.s., version of the analogy between BST and branching random walks, first given by Devroye \cite{Dev1}.

\section{Convergence of martingales}\label{marti}

In this section are given the main results about the asymptotic behaviors of  the Yule and BST martingales. The martingale connection (Proposition \ref{BHJ}) allows to express the links between the limits.

\subsection{Additive martingales}

Theorem \ref{cvu} gives an answer to a natural question asked in \cite{Jab2} 
about the domain {\it in the complex plane} where the BST martingale is $L^1-$convergent and uniformly convergent. Theorem \ref{optim} gives the optimal $L^1$ domain on $`R$.\par

\begin{theo}
\label{cvu}
For $1 < q <2$, let ${\cal V}_q :=\{ z: \sup_t `E |M(t, z)|^q < \infty\}$. 
Then
 ${\cal V}_q =\{ z: \ f(z, q) > 0\}$ with
\ben
\label{deff}
f(z, q) := 1 + q(2\Re z -1) -2|z|^q\,.
\een 
If we denote ${\cal V} := \cup_{1< q < 2} {\cal V}_q$, we have :
\begin{itemize}
\item[a)] As $t \rightarrow \infty$, $\{ M(t, z)\}$ converges, a.s. and in $L^1$, uniformly on every compact $C$ of ${\cal V}$. 

\item[b)]
As $n \rightarrow \infty$, $\{ `M_n (z)\}$ converges, a.s. and in $L^1$, uniformly on every compact $C$ of ${\cal V}$. 
\end{itemize}
\end{theo}

\proof a) is proved in \cite{JDB92} Theorem 6 (see also \cite{BertRou2}). \par
b) We will prove 
\begin{equation}
\label{25}
\lim_N \sup_{n \geq N} `E \sup_{z \in C} |`M_n (z) - `M_N (z)| = 0\,,
\end{equation}
which implies the uniform $L^1$ convergence and, 
since $(\sup_{z\in C} |`M_n (z) - `M_N (z)|)_{n \geq N}$ is a submartingale, 
this will imply also 
the a.s. uniform convergence\footnote{For the uniform a.s. convergence, 
it is possible to give a proof directly from \cite{JDB92}}. 
From  
the martingale connection (Proposition \ref{BHJ}), we have
$$`M_n (z) - `M_N (z) = `E [M(\tau_n ,z) - M(\tau_N , z) 
| {\cal F}_{(n)}]$$
so that taking supremum and expectation we get
$$`E \sup_{z \in C}|`M_n (z) - `M_N (z)|   \leq  `E \left( \sup_{z \in C} |M(\tau_n ,z)- M(\tau_N , z)|\right)\,.$$
Taking again the supremum in $n$ we get
\begin{eqnarray}
 \sup_{n \geq N} `E \sup_{z \in C}|`M_n (z) - `M_N (z)|   \leq  
 `E \sup_{n \geq N}\left( \sup_{z \in C} |M(\tau_n ,z)- M(\tau_N , z)|\right) \leq `E \Delta_n\,,
\end{eqnarray}
where we have set $\Delta_n := \sup_{T \geq \tau_n}\left( \sup_{z \in C} |M(T ,z)- M(\tau_n , z)|\right)$. 
Since $M(t,z)$ converges a.s. uniformly, we have 
a.s. $\lim_n \Delta_n = 0$. Moreover, by the triangle inequality $\Delta_n \leq 2  \Delta_0$, 
and by the proof of Proposition 1 in \cite{BertRou2}, $\Delta_0$ is integrable. 
The dominated convergence theorem gives $\lim_n `E \Delta_n = 0$ and (\ref{25}) holds, 
which ends the proof of Theorem \ref{cvu}.
\QED

\begin{rem}\label{BDG}  As usual the $L^1$ convergence in $a)$ of the above theorem comes from a $L^q$ bound (for some $1<q\leq 2$); more precisely, following the steps in \cite{Bertoin03} section 2.4, the quantity 
 $$\beta_t (\lambda) := \left( M(t, z) - 1\right) e^{t(2z-1)}$$
satisfies
\begin{equation}
\label{32}
E \mid\beta_{t} (z)\mid^q \leq e^{ t q(2\Re z -1)} 
\int_0 ^{t} \exp \left(-s f(z,q)\right) \,ds ~~~~\textrm{ for }1 < q  \leq 2.
\end{equation} 


\end{rem}

\begin{theo}
\label{bigjab} Let us assume $z \in (z_c ^- , z_c ^+)$.
\begin{itemize}
\item[1)] We have the 
{\bf limit martingale connection} :
\ben
\label{bhj}
\hbox{a.s.}\ \ \ M  (\infty, z) = \frac{\xi^{2z-1}}{\Gamma(2z)}\ `M _\infty (z)\,,
\een
where the exponential variable $\xi$ is defined in (\ref{limetnt}).

\item[2)] We have the following two splitting formulas:
\begin{itemize}
\item[a)] for the Yule process, 
\ben
\label{split}
M(\infty, z) = ze^{(1-2z)\tau_1} \left(M_0 (\infty, z) + M_1 (\infty, z)\right)
\ \ \ a.s.\een 
where $M_0 (\infty, z)$ and $M_1 (\infty, z)$ are independent, distributed as $M (\infty, z)$ and independent of $\tau_1$. 
\item[b)] for the BST,
\ben
\label{masterm}
`M _{\infty} (z) = z \left(U^{2z-1} `M _{\infty, (0)} (z) + (1-U)^{2z-1} `M _{\infty, (1)} (z)\right)
\een
where $U  \sim {\cal U}([0,1])$ is defined in (\ref{UUU}), $`M _{\infty ,(0)} (z), `M _{\infty, (1)} (z)$ are independent (and independent of $U$) 
and distributed as $`M _{\infty} (z)$.
\end{itemize}
\end{itemize}
\end{theo}
\proof 
1) is a consequence of (\ref{taun}) and  the martingale connection (\ref{mc}).       

2) a) For $t>\tau_1$ we have the decomposition
\ben
\label{firstsplit}
M(t, z) = z e^{(1-2z)\tau_1} \left[M^{(0)} (t-\tau_1 , z) + M^{(1)} (t-\tau_1 ,z) \right]
\een
where for $i= 0,1$ 
$$M^{(i)}(s, z) = \sum_{u \in \partial \BBt_s^{(i)}} z^{|u|} e^{s(1-2z)}\,,$$
and $\BBt^{(i)}$ is defined in Section \ref{bis}.

b) Take $t = \tau_n$ in (\ref{firstsplit}), condition on the first splitting time $\tau_1$, apply the branching property, let $n \rightarrow \infty$  
and apply  the limit martingale connection (\ref{bhj}) to get
\ben
\frac{\xi^{2z-1}}{\Gamma(2z)}`M_\infty  (z)= z e^{(1-2z)\tau_1}\left(\frac{\xi_0^{2z-1}}
{\Gamma(2z)}`M_{\infty,(0)}  (z) + \frac{\xi_1^{2z-1}}{\Gamma(2z)}`M_{\infty, (1)}  (z)\right)
\een
where $\xi_0$ and $\xi_1$ come from section \ref{bis},
which yields b) 
with the help of (\ref{xiutau}).
\QED

\noindent The following theorem gives the behavior in the remaining cases
\begin{theo}\label{optim}
For $z \in (0, \infty) \setminus 
(z^- _c , z^+ _c)$, then a.s. $\lim_t M(t , z) = 0$ and $\lim_n `M_n (z) = 0$.
\end{theo}
\proof
The continuous time result is in \cite{JDB92} (see also \cite{BertRou2}); it remains to use again the martingale connection (\ref{mc}).
\QED

\subsection{Derivative martingales}
From the above section, we deduce that 
the derivatives
\begin{equation}
M' (t ,z) := \frac{d}{dz} M(t,z), \ \ `M'_n (z) := \frac{d}{dz} `M_n (z)
\end{equation}
are  martingales which are no longer positive. They are called the derivative martingales. Their 
behaviors  are ruled by the following theorem. 
\begin{theo}
\label{theoderiv}
\begin{itemize}
\item[1)]  For $z \in (z^- _c , z^+ _c)$, the martingales 
$(M' (t ,z), t \geq 0)$ and $(`M'_n (z), n \geq 0)$ are convergent a.s..  
Let us call $M'(\infty, z)$ and $`M'_\infty (z)$ their limits. 
\begin{itemize}
\item[a)] We have the (derivative martingale) connection:
\ben
\label{relderiv}
M' (\infty, z) = \frac{\xi^{2z-1}}{\Gamma(2z)} \left(`M'_\infty (z)  + 
2 \left(\log \xi - \frac{\Gamma'(2z)}{\Gamma(2z)} \right)`M _\infty (z)\right)\, \hskip 1cm a.s.
\een
where $\xi \sim {\bf Exp}(1)$ is defined in (\ref{limetnt}) 
and is independent of $`M _\infty (z)$ and $`M'_\infty (z)$. 
\item[b)] We have the splitting formula:
\ben
\label{master'}
`M'_\infty (z) &=&  zU^{2z - 1} `M'_{\infty ,(0)} (z) + z(1-U)^{2z-1}`M'_{\infty ,(1)} (z) \cr
&+&  2z  \left(U^{2z-1} \log U\right) `M_{\infty,(0)} (z) + 2z \left((1-U)^{2z-1} \log (1-U)\right)`M_{\infty,(1)} (z)\cr
&+&
 z^{-1} `M_\infty (z) 
\een
where $U \sim {\cal U}([0,1])$ is defined in (\ref{UUU}), and the r.v. $`M'_{\infty ,(0)} (z)$ and $`M'_{\infty, (1)} (z)$ are independent 
(and independent of $U$) 
and distributed as $`M'_{\infty} (z)$.
\end{itemize}
\item[2) a)] The  martingales $(M' (t ,z^-_c), t \geq 0)$ and $(`M'_n (z^-_c), n \geq 0)$  (resp. $(M' (t ,z^+_c), t \geq 0)$ and $(`M'_n (z^+_c), n \geq 0)$)
are convergent a.s..
 Their limits denoted by $M' (\infty ,z^-_c)$ and $`M'_\infty (z^-_c)$ (resp. $M' (\infty ,z^+_c)$ 
 and $`M'_\infty (z^+_c)$) are positive (resp. negative) and satisfy
\ben
`E\big(M'(\infty,z_c^-)\big) = `E(`M_\infty(z_c^-))= +\infty,\\
`E\big(M'(\infty,z_c^+)\big) = `E(`M_\infty(z_c^+))= -\infty.
\een
\item[b)] $M' (\infty ,z^\pm_c)$ and $`M'_\infty (z^\pm_c)$ satisfy equations similar to 
(\ref{bhj}), (\ref{split}) and (\ref{masterm}):
\ben
\label{bhjc}
M'(\infty, z^\pm _c) &=&\frac{\xi^{2z^\pm _c -1}}{\Gamma(2z^\pm _c)}\ {\cal M}' _\infty (z^\pm _c)\\ 
\label{split'}
M'(\infty, z^\pm _c) &=& z^\pm _c e^{(1-2z^\pm _c)\tau_1} \left(M'_0 (\infty, z^\pm _c) + 
M'_1 (\infty, z^\pm _c)\right)
\\
\label{qsz}
`M'_{\infty} (z^\pm _c) &=& z^\pm _c \Big(U^{2z^\pm _c -1} `M'_{\infty, (0)}(z^\pm _c) + (1-U)^{2z^\pm _c -1} `M'_{\infty ,(1)} (z^\pm _c)\Big)\ \ \hbox{a.s.}\,.
\een
\end{itemize}
\end{theo}
\medskip

\proof 
1)  For $z \in (z_c ^- , z_c ^+)$ the a.s. convergence of $M'(t,z)$ is a consequence of the uniform convergence of $M(t,z)$ (by Theorem \ref{cvu}) and analyticity. 
Taking derivatives in  the martingale connection 
(\ref{mc})
 gives
\ben
\label{equaderiva}
M'  (\tau_n , z) = \left[\frac{C'_n (z)}{C_n(z)} - 2\tau_n\right]
 \textrm{\cls{C}}_n (z) `M_n (z)
+
 \textrm{\cls{C}}_n (z) `M'_n (z)\,.
\een
Using (\ref{etn}) again and
$$\frac{C'_n (z)}{C_n (z)} = \sum_{j=0}^{n-1}\frac{2}{j+2z} \ \ , \ \
\frac{\Gamma' (x)}{\Gamma(x)}
= \lim_n \Big(\log n - \sum_{j= 0}^{n-1} \frac{1}{x+j}\Big)\,,$$
we get
\ben
\label{cprimen}
\hbox{a.s.} \ \lim_n \left[\frac{C'_n (z)}{C_n(z)} - 2\tau_n \right] &=& 2 \left[-\frac{\Gamma'(2z)}{\Gamma(2z)} + \log \xi\right]\,.
\een
We conclude that $`M_n' (z)$ converges and that ${\cal M}'_{\infty}(z)$
satisfies (\ref{relderiv}) which proves a). 

To prove b),
we differentiate (\ref{firstsplit}) with respect to $z$ 
$$
M'(t,z) = (z^{-1} - 2\tau_1) M(t,z) +  z e^{(1-2z)\tau_1} \left[M^{(0)'} (t-\tau_1 , z) + M^{(1)'} (t-\tau_1 ,z) \right]\,, 
$$
and we use the same technique as above:
take $t=\tau_n$, let $n \rightarrow \infty$, apply (\ref{relderiv}) and its analogs with $(M'^{(i)}, {\cal M}^{(i)}, {\cal M}'^{(i)}, \xi_i)_{i =0,1}$ instead of $(M', {\cal M}, {\cal M}', \xi)$, and use (\ref{xiutau}).
\smallskip

\noindent 2) For $z = z^\pm  _c$, the a.s. convergence of the martingales $M'(t,z)$ and the signs of the limits are proved in \cite{BertRou2}, and so is the relation
$$`E M'(\infty , z_c^-)= - `EM'(\infty, z_c^+) = \infty\,.$$ 

Relation (\ref{bhjc}) is a consequence of (\ref{equaderiva}) and (\ref{cprimen}), since $`M_\infty   (z_c^\pm) = 0$.

Formula (\ref{qsz}) of 2)  is straightforward from (\ref{master'}) since $`M_\infty (z_c ^\pm ) =0$. Formula (\ref{bhjc}) is (\ref{relderiv}) for $z = z_c^\pm$.  
\QED
An easy but interesting consequence of  (\ref{master'}) is obtained 
in the following corollary, just taking  $z =1$ in (\ref{relderiv}) and (\ref{master'}) (remember that ${\cal M}_n (1)\equiv 1$). The 
distributional (weaker) version  of (\ref{quicksort}) below is the subject of a broad literature (see for instance Fill, Janson, Devroye, Neininger, R\"osler, R\"uschendorf \cite{FillJanson01,FillJanson02,DFN, NR, Ros,ROE}) and some properties of the distribution of $`M'_{\infty} (1)$ remain unknown.
\begin{cor}\label{cor11}We have
\ben
M'  (\infty , 1) = \xi \left(`M'_\infty (1) + 2 \left(\log \xi +  \gamma - 1\right)\right)\ \ 
\hbox{a.s.}\,,
\een
where $\gamma$ is the Euler constant, and $`M'_\infty (1)$ 
satisfies the a.s. version of the {\bf Quicksort} equation:
\ben
\label{quicksort}
`M'_\infty (1) = U `M'_{\infty ,(0)} (1) + (1-U) `M'_{\infty ,(1)} (1)
+ 2 U \log U + 2 (1-U) \log (1-U) + 1\,,
\een
where as above,  
$`M'_{\infty ,(0)} (1)$ and $`M'_{\infty, (1)} (1)$ are independent 
(and independent of $U$), distributed as $`M'_{\infty} (1)$ and $U \sim {\cal U}([0,1])$.
\end{cor}

\section{Convergence of  profiles}
\label{prof}
\label{profile}
\subsection{Random measures and profiles}
Recall that the profile of the tree ${\cal T}_n$  is the sequence
$$U_k (n) = \#\{u \in \partial {\cal T}_n : |u| = k\}, \ \ k \geq 1
\,,$$ 
and that, according to (\ref{cc'}), for every $\epsilon > 0$, there exists a.s. $n_0$ such that for  $n \geq n_0$, 
$$U_k(n) = 0 ~~ \textrm{ for }~~  k \notin [(c'-\epsilon)\log n , (c+\epsilon)\log n]\,.$$
It means that the convenient scaling for $k$ is  $(\log n)^{-1}$.
We are interested in the asymptotic behavior of $U_k (n)$ for $ k\cong x \log n$ and $x$ fixed in $(c',c)$.
 It is well known that 
$$`E\big(U_k (n)\big) = \frac{2^k}{n!} S_n ^{(k)}$$
where $S_n ^{(k)}$ is the Stirling number of the first kind. By analysis of singularities, 
Hwang (\cite{Hwang}) got an asymptotic estimate; for any $\ell > 0$ as $n\rightarrow \infty$ and 
$k\rightarrow \infty$ such that $r = k/ \log n \leq \ell$ :
\ben
\label{hwang+}
`E\,U_k (n) = \frac{(2\log n)^k}{k!\, n\, \Gamma(r)} (1 + o(1))\,.
\een
We deduce easily that for any $\ell > 0$ :
\ben
\label{hwang++}
`E\,U_k (n) =\frac{n^{1 -\eta_2 (\frac{k}{\log n})}}{\Gamma(\frac{k}{\log n})\sqrt{2\pi k}} (1 + o(1))\,,
\een
where $o(1)$ is uniform for $k/ \log n \leq \ell$ and $\eta_2$ was defined in (\ref{defeta}).
 
Jabbour in \cite{Jab1} introduced the random measure counting the levels of leaves in ${\cal T}_n$
\be
r_n := \sum_{k} U_k (n) \delta_k
 = \sum_{u \in \partial{\cal T}_n} \delta_{|u|}\,.
\ee
He proved that 
for $x \in (2 , c)$ 
\ben
\label{rvd++}
\hbox{a.s.} \ \ \lim_{n \rightarrow \infty} \frac{1}{\log n}\, \log  r_{n} ( ]x\log n, \infty[) &=& 1 -\eta_2 (x)
\een
and that the same result holds for $x \in (c' , 2)$, replacing 
$]x\log n, \infty[$
by 
$]0 , x\log n[$
.

At the level of random variables, Jabbour \& al. proved in \cite{Jab2} that 
\begin{equation}
\label{cvprof}
\hbox{a.s.}\ \ \lim_n \sup_{k/\log n \in[1.2 ,\  2.8]}\left(\frac{U_k (n)}{`EU_k (n)} - `M_\infty \big(\frac{k}{2\log n}\big)\right) = 0\,.
\end{equation}
Since their approach relies on $L^2$ estimations of $`M_n (z)$ with $z = k/(2\log n)$,
they guessed that the range $[1.2 , \ 2.8]$ may be extended to 
$I := (2 - 2^{1/2} , 2 + 2^{1/2})= (0.585... \ , 3.414...)$ which corresponds to the maximal interval in $z$
of $L^2$ convergence.
 In the following  subsection we extend the validity of the above result to (compact subsets of) 
the entire 
 interval $(c' , c) = (0.373... \ , 4.311...)$. 

This type of result is very reminiscent of sharp large deviations  in branching random walks
(\cite{JDB79}, \cite{JDB92}, \cite{Uchi}). Actually, we use the embedding method 
and results on the Yule process, as a branching random walk in continuous time.
The random measure counting the levels of leaves in the Yule tree is
$$\rho_t  = \sum_{u \in \partial \BBt_t} \delta_{|u|}\,,$$
(recall (\ref{defrho})). 
With the notations of \cite{Uchi}, 
the exponential rate of growing is ruled by the function 
$$x \mapsto L^\star (x) :=\sup_\theta \theta x - L(\theta) = \eta_2 (x) - 1\,,$$where the function $L$ is defined in \eref{lth}. This allows to define three areas:

\noindent-- for $x \in ]c' , c[$ , $\eta_2 (x) < 1$, so there are in mean about $e^{(1 - \eta_2 (x))t}$ leaves at level $\simeq xt$.
Call  this interval $]c' , c[$ ``supercritical area''.\\
\noindent -- for $x \in [0, c'[ \cup ]c, \infty[$, $\eta_2 (x) > 1$, so there are in mean about $e^{(1 - \eta_2 (x))t}$ leaves
at level  $\simeq xt$. Call this set ``subcritical area''.\\
\noindent -- call the set $\{c' , c\}$ ``critical area'' .

More precisely (Theorem 1' p. 909 \cite{Uchi}), for $x$ in the supercritical area,
\ben
\label{uchi1}
\lim_{t \rightarrow \infty} \sqrt{t}\ e^{tL^{\star}(x)} \rho_t ([xt]) = 
\sqrt{\frac{(L^{\star})'' (x)}{2 \pi}} \ M(\infty, x/2)\ \ \hbox{a.s.}\,.
\een
 

It is now tempting to replace $t$ by $\tau_n$ and $\rho_t ([xt])$ by $\rho_{\tau_n} ([x\log n]) = r_n([x\log n])$.
To validate this, we need some uniformity in $x$ in (\ref{uchi1}). In \cite{JDB92}, Biggins obtained such a result.
 However, it was in the non-lattice case, so we give in the next subsection a complete proof.
\subsection{Main result}
The aim of this section is to prove the following result
\begin{theo}\label{theopro} Almost surely, for any compact subset $K$ of $(c' , c)$
\begin{equation}
\lim_n \sup_{k : (k/\log n) \in K}\Blp\frac{U_k(n)}{\E\blp U_k(n)\brp}-\mathcal{M}_\infty\blp\frac{k}{2\log n} \brp \Brp=0\,.
\end{equation}
\end{theo}

\proof The following lemma, whose proof is postponed yields an asymptotic uniform behavior for $\rho_t (k)$.
\begin{lem}\label{promes}
Almost surely, for any compact $C$ of $(z^-_c , z^+_c)$, 
\ben
\label{limsup}
\lim_{t \rightarrow \infty} \sup_{k\geq 1 , z \in C}
z^{k} \sqrt{t} e^{t(1-2z)} \big[\rho_t (k) - M(\infty, z)e^{-t}\frac{(2t)^k}{k!}\big] = 0\,.
\een
\end{lem}
\medskip

Let $C$ be a compact subset of $(z_c^-,z_c^+)$.
From Lemma \ref{promes}, we know that
\be
\rho_t (k) = M(\infty, z)e^{-t}\frac{(2t)^k}{k!} + o(1)z^{-k} t^{-1/2}e^{-t(1-2z)}\,.
\ee
Recall that $o(1)$ is uniform in $k$ and in $z\in C$. If ${\cal P}^{(\lambda)}$ stands for the Poisson law with parameter  $\lambda$, notice that a ${\cal P}^{(2t)}$ appears in the previous expression.
Using a change of probability from ${\cal P}^{(2t)}$ to ${\cal P}^{(2tz)}$, we get
\be
\rho_t (k) = z^{-k} t^{-1/2}e^{-t(1-2z)}\Big[ t^{1/2}M(\infty, z) {\cal P}^{(2tz)} (k) + o(1)\Big]\,.
\ee
Using the  local limit theorem \cite{Petr}, we have
\be
\lim_{\lambda \rightarrow \infty} \sup_k \Big|\sqrt{2 \pi \lambda} \ 
{\cal P}^{(\lambda)} (k) - \exp \Big(- \frac{(k - \lambda)^2}{2\lambda}\Big) \Big| =0\,. 
\ee
Now, we set $\lambda=2tz$ with $z \in C$ which yields
\be
\lim_{t \rightarrow \infty} \sup_{z\in C}\sup_k \Big|\sqrt{4 \pi tz} \ {\cal P}^{(2tz)} (k) - \exp \Big(- \frac{(k - 2tz)^2}{4tz}\Big) \Big| =0.
\ee
Hence,
\begin{eqnarray}
\label{r=a}
\rho_t (k) = A_t(k, z)\Big[\exp \Big(- \frac{(k - 2tz)^2}{4tz}\Big) M(\infty, z)    +\Blp(4\pi z)^{1/2}+  M(\infty, z)\Brp o(1)\Big],
\end{eqnarray}
with
\begin{eqnarray*}
A_t(k,z) := \frac{e^{-t(1-2z)}}{z^k(4\pi tz)^{1/2}}\,.
\end{eqnarray*}
Remembering that 
$U_k(n)= \rho_{\tau_n} (k)$, we take $t=\tau_n$ and $\displaystyle{z=\frac k{2\log n}}$ in (\ref{r=a}).
Using (\ref{etn}) again and the estimate (\ref{hwang++}), 
we get
\begin{eqnarray*}
\frac{A_{\tau_n} (k, z)}{[`E U_k (n)]\ \xi^{1-2z} \Gamma (2z)} = 1 + o(1)\ \ , \ \ 
\exp \Big(- \frac{(k - 2\tau_n z)^2}{4\tau_n z}\Big) = 1 + o(1)\,.
\end{eqnarray*}
Now we apply the limit martingale connection (\ref{bhj}) and notice that
$$\sup_{z \in C} \Blp(4\pi z)^{1/2}+  M(\infty, z)\Brp < \infty$$
and we conclude 
$$U_k (n) = [`E U_k (n)] `M_\infty(z) (1 + o(1))\,,$$
with $z = k /(2 \log n)$ and $o(1)$ uniform in $z \in C$. \QED
\subsection{Proof of Lemma \ref{promes}}
We use the following lemma, which is the continuous time version of Lemma 5 
in \cite{JDB92}. Its proof can be managed with the same arguments, replacing Lemma 6 there, by Remark \ref{BDG}. We omit the details.
\begin{lem}\label{lemes}
\label{w-w}
For any $z_0 \in (z^-_c , z^+_c)$  there exists $r>0$ for which  $z^-_c < z_0 -r < z_0 +r < z^+_c$
and such that a.s. 
\begin{equation}
\lim_{t \rightarrow \infty}\sup_{z \in [z_0 -r , z_0 +r]} \int_{-\pi}^{\pi}\sqrt t \mid M(t, ze^{i\eta} ) - M(\infty, z)\mid 
e^{-2tz(1 - \cos\eta)} d \eta = 0\,. 
\end{equation} 
\end{lem}

\medskip
Write
$$M(t,z) = e^{t(1-2z)} \sum_k \rho_t (k) z^k\,,
$$
and the Fourier inversion formula yields
$$
\rho_t (k) = \frac{e^{-t(1-2z)}z^{-k}}{2\pi}\int_{-\pi}^{\pi} M(t, ze^{i\eta}) e^{-2tz(1-e^{i\eta})}e^{-ik\eta}d\eta
$$ 
and, owing to Lemma \ref{lemes}
$$2\pi \rho_t (k) e^{t(1-2z)}z^k \sqrt{t} = M(\infty, z) \sqrt{t} \int_{-\pi}^\pi e^{-2zt (1- e^{i\eta})} e^{-ik\eta} d\eta
+ o(1)$$
with $o(1)$ uniform in $k$ and in $z$ in any compact subset of  $(z^-_c , z^+_c)$.
Now, from the Cauchy formula we get that 
$$
\int_{-\pi}^\pi e^{-2zt(1 - e^{i\eta})} e^{-ik\eta} d\eta = 2\pi e^{-2zt} \frac{(2zt)^k}{k!}\,, 
$$
yielding (\ref{limsup}), which ends the proof.
\QED
\section{Tagged branches and biased trees}
\label{tag}

We now introduce the tilting or biasing method in the setting of the Yule and BST processes.

This procedure consists  in marking at random 
a special ``ray'' or  branch of the Yule tree.
The special ray of the BST is the spine of the marked Yule tree observed in the splitting times. One then applies to this special ray some  evolution rule, 
different from the other branches. The result of this change of probability is that the 
whole tree owns a different behavior. This method is usual and fruitful in modern developments on 
branching processes, and also in the study of fragmentation processes (\cite{BertRou2,BK1,Lyons1}).
The introduction of the tilting method in the setting of BST 
provides new tools to study some characteristics of the BST. 
\par
The martingales studied above are the right change of probability to pass from a tilted model to the non-tilted model: they appear as  Radon-Nikodym derivatives. The parameter $z$, present in the martingales $(M(t,z))_{t\geq 0}$ and $({\cal M}_n(z))_{n\geq 0}$, allows to tune the growing of the special ray, changing in a visible way the shape of the (Yule or BST) tree.

\subsection{Tilted fragmentation and biased Yule tree}
\label{TF}

First at all, let us enlarge the probability space of the fragmentation process introduced in Section \ref{YT}. 
 Let us denote by ${\cal F}_t$  the $\sigma$-algebra of the interval 
fragmentation process $F(.)$ up to time $t$ and $V$ be a ${\cal U} ([0,1])$ r.v. independent of the filtration $\left({\cal F}_t\right)_{t\geq 0}$.\par

Since $`P( V \in \{ k2^{-j}, 0 \leq k \leq 2^j, j\in `N , k\in `N \}) = 0$, 
we may define $`P$-a.s. for every $t$ a unique $S(t) \in \mathbb{U}$ 
such that $I_{S(t)}$ is an interval of $ F(t)$ and $V \in I_{S(t)}$ . 
In other words, $S(t)$ is the element of $\mathbb{U}$ encoding
 the fragment containing $V$, its depth is $s(t) := |S(t)|$, the length of  $I_{S(t)}$ is 
 $2^{-s(t)}$ and
\ben
\label{sizeb}
`P (S(t) = u \ | \ {\cal F}_t) = 2^{-|u|}\, \, , \  u \in \partial\BBt_t
\een 
(it is equivalent to  choose a fragment at random with probability equal to its length,
it is the classical size-biasing setting).

Now we build the process $(\widetilde\BBt_t)_{t\geq 0}$ of marked binary Yule trees associated with the pair $(F(.),S(.))$. 
The only change with Section \ref{YT} is the role played by the random variable $V$ (missing in Section \ref{YT}). 
During the construction of the Yule tree, at any given time $t$, each leaf in $\BBt_t$ corresponds to an interval 
in the fragmentation $F(t)$. For every $t$ we mark the leaf $S(t)$ of $\BBt_t$  
that corresponds to the interval $I_{S(t)}$ that contains $V$. We obtain a marked tree called $(\widetilde\BBt_t)_{t\geq 0}$. Thus, the set of nodes marked during $[0,t]$ are the prefixes of $S(t)$. 
We call spine the process $S(.)$.

In fact,  given $\widetilde\BBt_t$, one can recover $(F(t),S(t) )$. Moreover, with the whole process $(\widetilde\BBt_t)_{t\geq 0}$  one can a.s. recover $V$: 
\[V=\bigcap_{t\geq 0} I_{S(t)}.\]

As a consequence of the general theory of homogeneous fragmentations (see Bertoin \cite{Bertoin01}) 
or by a direct computation, we see that $(s(t) , t \geq 0)$ is an homogeneous Poisson process with parameter $1$. In particular, if
\ben
\label{expcont}
{\cal E}(t,z) :=  (2z)^{s(t)
}
e^{t(1-2z)}
\een
then $`E \ {\cal E}(t,z) = 1$. Conditionally on 
$\widehat{\cal F}_r = {\cal F}_r \vee \sigma (S(r) , s \leq r)$, 
the restriction of the fragmentation $F(.+r)$ to the interval $I_{S(r)}$ is 
distributed as a rescaling of $F(.)$ by a factor $2^{-s(r)}$, 
which entails that $\big( {\cal E}(t,z),\widehat{\cal F}_t\big)_{t\geq 0}$ is a martingale.
By  the size biasing scheme (\ref{sizeb}) and the definition (\ref{BHM}) we get
\ben
\label{sbcont}
M(t,z) = `E\left[{\cal E}(t,z) \ | \ {\cal F}_t\right].
\een
Hence, the Yule martingale appears to be a projection of the martingale ${\cal E}$ (which is a spine-measurable function) on the $\sigma$-algebra containing only the underlying binary tree.

Coming  back to the discrete time, 
set {\bf Spine}$_n:= 
S(\tau_n)$ and $s_n := |\hbox{\bf Spine}_n |$. Notice that the underlying unmarked tree $\BBt_{\tau_n}$ is ${\cal F}_{(n)}$-measurable.\par
 Applying  (\ref{sizeb}) at the $(
{\cal F}_t , t \geq 0)$ stopping time
$\tau_n$, we get for every leaf $u \in \partial{\cal T}_n$ (and $k \geq 1$) :
\ben
\label{loisn}
`P(\hbox{\bf Spine}_n=u\ | \ {\cal F}_{(n)})& =& 2^{-|u|} ,\\ 
\nonumber `P(s_n=k\ | \ {\cal F}_{(n)}
)&=& U_k(n) 2^{-|k|}\ .
\een
Thus, for fixed $n$, to draw at random the marked tree $\tilde{\BBt}_{\tau_n}$, one may choose at first a binary tree ${\cal T}_n$, and then pick the marked leaf according to the conditional distribution \eref{loisn}. 
Let $\widehat{\cal F}_{(0)}$ be the trivial $\sigma$-algebra, and for $n \geq 1$ let $\widehat{\cal F}_{(n)}$ be the  $\sigma$-algebra 
 obtained from ${\cal F}_{(n)}$ by adjunction of 
$S(\tau_1), ... , S(\tau_n)$. Let us consider
${\cal E}_n (z) := `E\Big[{\cal E}( \tau_n ,z )\ | \ \widehat{\cal F}_{(n)}\Big] $
 (with $ {\cal E}_0 (z) := 1$). From Lemma \ref{BNT} a) we have $`E(e^{\tau_n (1-2z)})= C_n (z)^{-1}$ hence
\ben
\label{expdiscr}
{\cal E}_n (z) = (2z)^{\displaystyle s_n} \, C_n (z)^{-1}.
\een
From the martingale property of ${\cal E}(t ,z )$ 
and the definition of ${\cal E}_n (z)$ we see that $\big({\cal E}_n (z) , 
\widehat{\cal F}_{(n)}\big)$ 
is a martingale. Like in (\ref{sbcont}), we get easily
\ben
\label{sbdiscr} 
`M_n (z) =  `E\left[{\cal E}_n(z) \ | \ {\cal F}_{(n)}
\right],
\een
so that  the martingales $M(t,z)$ and $`M_n(z)$ are obtained from the ``exponential martingales'' ${\cal E}(z,t)$ and ${\cal E}_n (z)$ by projection. 

Moreover the martingale connection (\ref{mc}) may be seen 
as the projection on $\widehat{\cal F}_{(n)}$ of the relation 
$${\cal E}(\tau_n , z) = \textrm{\cls{C}}_n (z){\cal E}_n (z)\,.$$
Note that one may also obtain ${\cal M}_n(z)$ as
\[{\cal M}_n(z)=`E(M(\tau_n,z)| {\cal F}_{(n)});\]
this is a kind of integration with respect to the time.
All these martingales are precisely the main tool to tilt  probabilities.
 In particular we define $`P^{(2z)}$ on $(\widehat{\cal F}_t , t \geq 0)$
 by
\ben
\label{tiltcont}
`P^{(2z)}_{|_{\widehat{\cal F}_t}} = {\cal E}(t,z) \ `P_{|_{\widehat{\cal F}_t}} \,,
\een

By projection on $({\cal F}_t, t \geq 0)$, (\ref{tiltcont}) yields
\ben
`P^{(2z)}_{|_{{\cal F}
_t}} = M(t,z) \ `P_{|_{{\cal F}_t}}\,.
\een

If $^d`P$ (resp. $^d`P^{(2z)}$) is the restriction of $`P$ (resp. $`P^{(2z)})$ to $\vee_n \widehat{\cal F}_{(n)}$, the
discrete versions of the above relations are
\ben
\label{tiltdiscr}
^d`P^{(2z)}_{|_{\widehat{\cal F}_{(n)}}} 
= {\cal E}_n(z) \ ^d`P_{|_{\widehat{\cal F}_{(n)}}}\ \ , \ \ 
 ^d`P^{(2z)}_{|_{{\cal F}_{(n)}
 }} 
= `M_n (z)  \ ^d`P_{|_{{\cal F}_{(n)}
 }} \ .
\een

It turns out that $`P^{(2z)}$ can be seen as a probability on marked Yule trees.
This is the object of the following subsection.

\subsection{A biased Yule tree}
\label{BYT}

Recall the construction of the Yule tree process $(\BBt_t)_{t\geq 0}$ given in Section \ref{YT}. Each leaf $u$ of the current Yule tree owns a  ${\bf  Exp} (1)$-distributed clock. At its death, $u$ becomes an internal node, and two leaves $u0$ and $u1$ appear (with new  ${\bf  Exp} (1)$, independent of the other ones). \par
Let us consider now a model of marked binary tree $(\BBt_t^{\star})_{t\geq 0}$ defined as follows.
\par
In $\BBt_t^{\star}$ there are now two kinds of nodes: marked and unmarked. We denote by $(v,m)$ 
the node $v$ if it is marked, and by $(v,\bar{m})$ the node $v$ if it is unmarked. 
At time 0,  $\BBt_0^{\star}=\{(\emptyset,m)\}$.

Each unmarked leaf owns a ${\bf  Exp} (1)$-distributed clock. The marked leaf owns a ${\bf  Exp} (2z)$-distributed 
clock. Now the evolution of the tilted Yule tree is as follows:\\
$\bullet$ when an unmarked leaf $u$ dies, $u$ becomes an unmarked internal node, and two unmarked leaves $(u0,\bar{m})$ and $(u1,\bar{m})$ appear.\\
$\bullet$ when the marked leaf $u$ dies, $u$ becomes a marked internal node. Two leaves $u0$ and $u1$ appear. We mark at random $u0$ or $u1$ (equally likely) and let the other one unmarked.\\
The marked nodes form a branch in the tree. The behavior of the marked branch depends on the value of $2z$. If $2z>1$ then, the growing of the marked branch is faster than the other ones, when $2z<1$, the growing of the marked branch is slower. 
The depth of the marked leaf follows a Poisson process of rate $2z$. Notice that we have already met this Poisson process in the proof of Theorem \ref{theopro}.\par

It turns out that under $`P^{(2z)}$, the process $(\tilde{\BBt_t})_{t\geq 0}$ has the same distribution as the process $(\BBt_t^{\star})_{t\geq 0}$ (consider the spine as the marked leaf). For the underlying branching random walk this construction is classical \cite{Ath99},\cite{BK1},\cite{KYP04},.... For the fragmentation it can be found in \cite{BertRou2}. \par
\medskip

\noindent{\bf Remark}
Denote by $S^{\star}(t)$ the marked leaf in $\BBt_t^{\star}$, and consider 
\[V^{\star}=\bigcap_{t\geq 0} I_{S^{\star}(t)}.\]
By symmetry of the splittings, the random variable $V^{\star}$ is ${\cal U}[0,1]$-distributed and independent of the process of the underlying unmarked trees derived from $(\BBt_t^{\star})_{t\geq 0}$. Note $F^{\star}(.)$ the fragmentation process that is associated with the unmarked process derived from $\BBt_t^{\star}$.

The a.s. bijection between $(F(.), S(.))$ and $(F(.), V)$ (under $`P^{(1)}$) explained 
in the beginning of Section \ref{TF} is also valid between $(F^{\star}(.), S^{\star}(.))$ 
and $(F^{\star}(.), V^{\star})$. 
Hence, the law on marked fragmentation $(F(.),S(.))$ under $`P^{(2z)}$ (defined by \eref{tiltcont}) 
is the law of $(F^{\star}(.), S^{\star}(.))$. It follows that, under $`P^{(2z)}$ 
one may also build the spine by choosing at first a uniform random variable $V$ 
and follow the fragment containing $V$. 
This is not true in general when using  the tilting method. Usually, at each splitting of the marked fragment $M$, one 
has to choose the new marked fragment among the children of $M$, 
according to a rule depending on the size of these fragments. 
It cannot be summed up by the drawing of a random variable $V$, once for all as in our case, where sizes are equal.   
\medskip

\noindent According to the representation by $(\BBt_t^{\star})_{t\geq 0}$, the Yule tree owns a natural decomposition according to the marked branch. Let $u$ be a node of the marked branch. 
One of the nodes $u0$ or $u1$ does not belong to this marked branch. Assume that it is $u0$. 
Then, (up to a change of the time origin),\\
$\bullet$ the subtree rooted in $u0$ is a copy of the untilted Yule tree;\\
$\bullet$  the subtree rooted in $u1$ is a copy of the tilted Yule tree.\\
We can also see this process as a branching process with immigration, as presented in \cite{Tavar} (see also \cite{Pit} chap. 10 and \cite{HopF}).

\subsection{A biased BST model}
\label{BYBS}
The tilted Yule tree can also been stopped at time $\tau_n$ of the creation time of the $n$th internal node. Let $\widetilde {\cal T}_n$ be the obtained marked binary search tree. The discrete evolution is as follows:\\
$\widetilde {\cal T}_n$ is a complete binary tree with $2n+1$ nodes, in which one leaf is marked and the $n$ other ones are unmarked. 
Knowing $\widetilde {\cal T}_n$, the marked tree $\widetilde {\cal T}_{n+1}$ is as follows:\\
we choose the marked leaf with probability $2z/(n+2z)$ and 
each unmarked one with probability $1/(n+2z)$.\\
$\bullet$ If the chosen leaf  $v$ is unmarked, then $v$ becomes an unmarked internal node and two unmarked leaves $v0$ and $v1$ are created.\\
$\bullet$ If the chosen leaf  $v$ is marked, $v$ becomes a marked internal node. Two leaves $v0$ and $v1$ appear. One marks at random $v0$ or $v1$ (equally likely) and let the other ones unmarked.\par
We note $\mathbb{Q}^{(2z)}$ for the law on the marked binary search tree process $(\widetilde {\cal T}_n)_n$ under this model of evolution.\par

Once again, the BST can be decomposed along the marked branch. 
The speed of growing of the marked branch depends on the value of $2z$. 
One may also interpret the size of the subtrees rooted on the tilted branch as tables in a Chinese restaurant 
(see Barbour \& al. \cite{ABT}, Pitman \cite{Pit}), and obtain, like this, new explanations of the behavior of the size of the subtrees rooted on the marked branch. \par

As in the previous subsection, we denote by $(v,m)$ a marked node and 
$(v,\bar{m})$
an unmarked node. The dynamics we described above yields the following conditional probabilities:

if $(v \bar m) \in \partial\widetilde{\cal T}_n$, then
$$\mathbb{Q}^{(2z)} ({\bf Spine}_{n+1} = {\bf Spine}_{n} , \widetilde{\cal T}_{n+1} = \widetilde{\cal T}_n \cup \{(v0, \bar{m}), (v1, \bar{m})\} |\widetilde{\cal T}_n ) = \frac{1}{n +2z}$$
If $(v,m) \in \partial\widetilde{\cal T}_n$, (i.e. ${\bf Spine}_{n} = v$), then
$$\mathbb{Q}^{(2z)} ({\bf Spine}_{n+1} = v0 , \widetilde{\cal T}_{n+1} = \widetilde{\cal T}_n \cup \{(v0, m), (v1, \bar{m})\} |\widetilde{\cal T}_n ) = \frac{1}{2}\frac{2z}{n +2z}$$ 
similarly,
$$\mathbb{Q}^{(2z)} ({\bf Spine}_{n+1} = v1 , \widetilde{\cal T}_{n+1} = \widetilde{\cal T}_n \cup \{(v0, \bar{m}), (v1, m)\} |\widetilde{\cal T}_n ) = \frac{1}{2}\frac{2z}{n +2z}\,.$$ 
Summing up, we have for any marked tree $\widetilde{t}_{n+1}$ with $n+1$ nodes that can be obtained from $\widetilde{\cal T}_n $ by one insertion
\ben
\mathbb{Q}^{(2z)}(\widetilde{\cal T}_{n+1} = \widetilde{t}_{n+1}|\widetilde{\cal T}_n ) = \frac{z^{s_{n+1} - s_n}}{n+2z}\een
and 
\[\mathbb{Q}^{(1)}(\widetilde{\cal T}_{n+1} = \widetilde{t}_{n+1}|\widetilde{\cal T}_n )= \frac{(1/2)^{s_{n+1} - s_n}}{n+1}.\]
Thus, by iterative construction,
\[\frac{\mathbb{Q}^{(2z)}}{\mathbb{Q}^{(1)}}\Big|_{{\hat{\cal F}_n}}=\prod_{j=0}^{n-1}  
\frac{(2z)^{s_{j+1} - s_j}(j+1)}{j+2z} = (2z)^{\displaystyle s_n} \, C_n (z)^{-1}={\cal E}_n (z)\,.\] 
Hence, $\mathbb{Q}^{(2z)}$ is absolutely continuous with respect to $\mathbb{Q}^{(1)}$, with the  Radon-Nikodym derivative announced in (\ref{tiltdiscr}). Since $\mathbb{Q}^{(1)}$ and $^d \BBp^{(1)}$ (the non-biased models) are identical, 
the law of $(\widetilde{\cal T}_n )_n$ under $\mathbb{Q}^{(2z)}$ is $^d \BBp^{(2z)}$. \par
One finds an analogous result (in another context) and its proof in Lemma 1 and 2 of \cite{BK1}.

\subsection{Spine evolution}
Thanks to the previous subsections, it 
appears that under $^d`P^{(2z)}$
\ben
\label{sber0}
s_n = 1 + \sum_1 ^{n-1} \epsilon_k
\een
where $(\epsilon_k)_{k\geq 1}$ are independent and for every $k \geq 1$, $\epsilon_k$ is a Bernoulli random variable with parameter $\frac{2z}{k+ 2z}$; (we use the notation $\epsilon_k \sim {\bf Be}(\frac{2z}{k+ 2z})$). 

\begin{prop}
\label{spineprop}
For any parameter $z > 0$,
\begin{itemize}
\item[1)] (strong law)
\ben
\label{SLLN}
\lim \frac{s_n}{\log n} = 2z\,,\  \ \ \ ^d`P^{(2z)}-\hbox{a.s.}.
\een
\item[2)] (central limit theorem)
The distribution of $\,\,\displaystyle\frac{s_n -2z\log n}{\sqrt{2z \log n}}\,\,$  under 
$^d`P^{(2z)}$ converges to a standard normal distribution ${\cal N}(0,1)$.
\item[3)] (large deviations)
\label{PGD}
The family of distributions of $(s_n , \ n > 0)$ under $^d`P^{(2z)}$ satisfies the large deviation principle on $[0, \infty)$ with speed $\log n$ 
and rate function $\eta_{2z}$ where the function  $\eta_\lambda$ is defined in (\ref{defeta}).
\end{itemize}
\end{prop}

\proof 
1) and 2) are consequences of known results on sums of independent r.v. (see \cite{Petr}). Notice also
 that 
$s_n - `E^{(2z)}(s_n)$ is a martingale.
 
 3) is a consequence of G\"artner-Ellis theorem.\QED
 
Once again, this proposition shows that under the biased model, the BST evolves rather differently that under the usual model. For example, the marked leaf depth is about $2z\log n$. So, for $z>z_c^+$, the marked leaf is higher that the height of the non-biased BST.

\subsection{Depth of insertion}

In introducing the BST model, we defined
 the sequence $(D_n , n \geq 0)$ 
as the successive inserted nodes and $d_n = |D_n|$ (see (\ref{transit})). 
In continuous time, we set $\eta (t) = \inf \{ s > t : \ \BBt_s \not= \BBt_t \}$ for 
the first time of growing after $t$, and $D(t)= \BBt_{\eta (t)}\setminus \BBt_t$ for the node of the coming insertion.

Let us stress on  the difference between the spine processes $(s_n , n \geq 0)$  and $(s(t) , t \geq 0)$ and the  insertion processes $(d_n , n \geq 0)$ and $(d(t) , t \geq 0)$.

The (marginal) distribution of $d_n$ is given in Jabbour \cite{Jab1} (see also Mahmoud \cite{Mah})
\ben
\label{socalledjab}
`E z^{d_n} = \frac{C_n (z)}{n+1} = \frac{(2z) (2z + 1) \cdots (2z + n-1)}{(n+1)!}
\een
so that
$$d_n \buildrel{law}\over{=} 1 + \sum_1^{n-1} \varepsilon_k\,,$$
where  $(\varepsilon_k)_{k\geq 1}$ are independent and for every $k \geq 1$, $\varepsilon_k \sim {\bf Be}(\frac{2}{k+ 2})$ .

\begin{prop} The following convergences hold:\\
$(i)$ $(d_n)$ satisfies a law of large numbers:
\ben
\frac{d_n}{2\log n} \buildrel{P}\over{\longrightarrow} 1;
\een
$(ii)$ it satisfies a central limit theorem:
\ben
\frac{d_n - 2 \log n}{\sqrt {2 \log n}} \buildrel{{law}}\over{\Longrightarrow} {\cal N}(0,1)\,.
\een
$(iii)$ We have
\ben
\label{limsupinf}
\hbox{a.s.}\ \liminf_n \frac{d_n}{2 \log n} = \frac{c'}{2} = z^-\  \ , \ \  
\hbox{a.s.}\ \limsup_n \frac{d_n}{2 \log n} = \frac{c}{2} = z^+\,.
\een 
\end{prop}
Note that $(iii)$ of course, implies that $(ii)$ is not an almost sure convergence.\\ 
\proof
The arguments to prove $(i)$ and $(ii)$ are classical; $(iii)$ is a consequence of  (\ref{cc'}).
\QED
For the Yule tree, we did not find the distribution of $d(t)$ in the literature. Let us give the joint distribution of $(N_t , d(t))$ (for $t$ fixed).

Since $\{ N_t =n +1 \} = \{ \tau_{n} \leq t < \tau_{n+1} \}$, we have
$`E ( z^{d(t)} s^{N_t}) = \sum_0 ^\infty (`Ez^{d_n}) `P(N_t =n+1)s^{n+1}  $.
Since the distribution of $N_t$ is geometric of parameter $e^{-t}$, and owing to (\ref{socalledjab}) we get
\ben
`E ( z^{d(t)}s^{N_t}) &=& \sum_0 ^\infty \frac{(2z) (2z + 1) \cdots (2z + n-1)}{(n+1)!} e^{-t} (1 - e^{-t})^{n}s^{n+1}\\
&=& \frac{\big(1 - s(1 - e^{-t})\big)^{1 -2z} - 1}{(e^t - 1)( 2z - 1)}\,.
\een
Taking $s=1$, we get the marginal of $d(t)$
$$`E  z^{d(t)} = \frac{e^{t(2z-1)} - 1}{(e^t - 1)( 2z - 1)}\,.$$ 
Transforming these generating functions into Fourier transforms, it is now easy to conclude that
\begin{prop}
As $t\rightarrow \infty$,
$$\left(N_t e^{-t}, \frac{d(t) - 2t}{\sqrt {2t}}\right) 
\buildrel{{law}}\over{\Longrightarrow} (\xi, G)$$
where $\xi$ is defined in (\ref{etn}) and $G$ is ${\cal N}(0,1)$ and independent of 
$\xi$.
\end{prop}
\bigskip

\noindent{\bf Remark}:  For the same reasons as in (\ref{limsupinf}), we have
$$
\hbox{a.s.}\ \liminf_t \frac{d(t)}{2t} = \frac{c'}{2} = z^-\  \ , \ \  
\hbox{a.s.}\ \limsup_t \frac{d(t)}{2t} = \frac{c}{2} = z^+\,.
$$
Under the change of probability $`P^{(2z)}$ (or using Kolmogorov equations) the distribution of $N_t$ is given by:
\ben
\label{loisousz}
`E^{(2z)} \phi^{N_t}&=& `E \big[(2z)^{s(t)} e^{t(1-2z)}\phi^{N_t}\big]\\
&=&  e^{t(1-2z)}\big[\phi + \sum_0 ^\infty (`E(2z)^{s_n}) `P(N_t =n+1)\phi^{n+1}\big ]\\
&=&  e^{t(1-2z)}\big[\phi e^{-t} + e^{-t} \sum_{n=2}^\infty (1-e^{-t})^{n-1} \prod_0^{n-2} 
\frac{j+2z}{j+1} \phi^n\big]\\
&=& \phi \left[\frac{e^{-t}}{1-\phi(1-e^{-t})}\right]^{2z}\,;
\een
where $\phi$ is any real in $[0,1]$. Hence, under $`P^{(2z)}$, the r.v. $N_t - 1$ is a negative binomial of order $2z$ and parameter $e^{-t}$.  
As $t\rightarrow \infty$, the $`P^{(2z)}$ distribution of 
$e^{-t}N_t$ converges to a  $\gamma(2z)$-distributed random variable.  
 Actually we have for every $z, t ,h$
\ben
`E[N_{t+h} | {\cal F}_t] &=& (N_t - 1) `E N_h + `E^{(2z)} N_h\\
&=& e^h (N_t - 1) + 1 + 2z (e^h - 1)\,.
\een
This implies that $(e^{-t}[N_t - 1 + 2z])_{t \geq 0})$ is a $`P^{(2z)}$ martingale. 
If $2z > 1$, it is positive hence convergent $`P^{(2z)}$ a.s. If
$2z < 1$, then $(e^{-t}N_t)_{t \geq 0})$ is a positive supermartingale, hence $`P^{(2z)}$ a.s. convergent.\QED
\bibliographystyle{plain} 
\bibliography{versionlongue}
\end{document}